\documentclass[reqno,onecolumn,12pt]{amsart}
\usepackage{amsfonts}
\usepackage{amssymb}
\usepackage{amsmath}
\usepackage{graphicx}
\usepackage{amscd,amsthm}

\newtheorem{theorem}{Theorem}
\theoremstyle{plain}

\newtheorem{corollary}{Corollary}

\numberwithin{equation}{section}

\input{tcilatex}

\begin{document}
\author{}
\title{}
\maketitle

\begin{center}
\pagestyle{myheadings} \thispagestyle{empty} 
\markboth{\bf Serkan Araci, Armen Bagdasaryan, Cenap \"{O}zel and H. M. Srivastava}
{\bf Some New Symmetric Identities for the $q$-Zeta Type Functions}

\textbf{\large Some New Symmetric Identities for the $q$-Zeta Type Functions}

\bigskip

\textbf{Serkan Araci$^{1,\ast}$, Armen Bagdasaryan$^{2}$, Cenap \"{O}zel$%
^{3} $ and H. M. Srivastava$^{4}$}

\bigskip

$^{1}$Atat\"{u}rk Street, TR-31290 Hatay, Turkey\\[0pt]

\textbf{E-Mail: mtsrkn@hotmail.com}\\[0pt]

\textbf{$^{\ast}$Corresponding Author}\\[2mm]

$^{2}$Institute for Control Sciences, Russian Academy of Sciences, 65
Profsoyuznaya Street, 117997 Moscow, Russian Federation\\[0pt]

\textbf{E-Mail: abagdasari@hotmail.com}\\[2mm]

$^{3}$Department of Mathematics, Abant Izzet Baysal University, TR-14268
Bolu, Turkey\\[0pt]

\textbf{E-Mail: cenap.ozel@gmail.com}\\[2mm]

$^{4}$Department of Mathematics and Statistics, University of Victoria,
Victoria, British Columbia V8W 3R4, Canada\\[0pt]

\textbf{E-Mail: harimsri@math.uvic.ca}\\[4mm]

\textbf{\large Abstract}
\end{center}

\begin{quotation}
The main object of this paper is to obtain several symmetric properties of
the $q$-Zeta type functions. As applications of these properties, we give
some new interesting identities for the modified $q$-Genocchi polynomials.
Finally, our applications are shown to lead to a number of interesting
results which we state in the present paper.
\end{quotation}

\begin{center}
\textbf{{\large R\'{e}sum\'{e}}}
\end{center}

\begin{quotation}
L'objet principal de cet article est d'obtenir plusieurs propri\'{e}t\'{e}s
sym\'{e}triques des fonctions de type $q$-Zeta. Comme les applications de
ces propri\'{e}t\'{e}s, nous donnons quelques nouvelles identit\'{e}s int%
\'{e}ressantes pour les polyn\^{o}mes $q$-Genocchi modifi\'{e}s. Enfin, nos
applications sont pr\'{e}sent\'{e}s \`{a} conduire \`{a} un certain nombre
de r\'{e}sultats int\'{e}ressants que nous disons dans le pr\'{e}sent
document.
\end{quotation}

\bigskip



\section{\textbf{Introduction}}

Throughout this paper, we use the following standard notations: 
\begin{equation*}
\mathbb{N}:=\left\{ 1,2,3,\cdots \right\} \qquad \text{ and }\qquad \mathbb{N%
}_{0}:=\left\{ 0,1,2,\cdots \right\} =\mathbb{N}\cup \left\{ 0\right\} \text{%
.}
\end{equation*}%
Also, as usual, $\mathbb{R}$ denotes the set of real numbers and $\mathbb{C}$
denotes the set of complex numbers.

The Genocchi polynomials $G_{n}(x)$ and the Genocchi numbers $G_{n}:=G_{n}(0)
$ are given by the following generating functions: 
\begin{equation}
\sum_{n=0}^{\infty }G_{n}(x)\frac{t^{n}}{n!}=\left( \frac{2t}{e^{t}+1}%
\right) e^{xt}\qquad \text{and}\qquad \sum_{n=0}^{\infty }G_{n}\;\frac{t^{n}%
}{n!}=\frac{2t}{e^{t}+1}\qquad (\left\vert t\right\vert <\pi )\text{,}
\label{e1}
\end{equation}%
respectively. In particular, the second generating function in \eqref{e1}
can be restated as follows: 
\begin{equation*}
e^{Gt+t}+e^{Gt}=2t
\end{equation*}%
by using the \textit{umbral }(symbolic) convention exhibited by $G^{n}:=G_{n}
$. By utilizing the Taylor-Maclaurin expansion, one finds that 
\begin{equation}
\left( G+1\right) ^{n}+G_{n}=\left\{ 
\begin{array}{ll}
2 & \qquad (n=1) \\ 
&  \\ 
0 & \qquad (\mathrm{otherwise})\text{.}%
\end{array}%
\right.   \label{e2}
\end{equation}%
It follows from (\ref{e2}) that (see, for details, \cite{SrivastavaChoi2012}%
) 
\begin{equation*}
G_{1}=1,\;\;G_{2}=-1,\;\;G_{3}=0,\;\;G_{4}=1,\;\;G_{5}=0,\;\;G_{6}=-3,\;%
\;G_{7}=0,\;\;G_{8}=17,\cdots 
\end{equation*}%
and (in general) 
\begin{equation*}
G_{2n+1}=0\qquad (n\in \mathbb{N})\text{.}
\end{equation*}

The history of the Genocchi polynomials $G_{n}(x)$ and the Genocchi numbers $%
G_{n}$ can be traced back to the Italian mathematician, Angelo Genocchi
(1817--1889). From Genocchi to the present time, the Genocchi polynomials
and the Genocchi numbers have been extensively studied in many different
contexts in such branches of Mathematics as, for instance, Elementary Number
Theory, Complex Analytic Number Theory, Homotopy Theory (especially stable
Homotopy groups of spheres), Differential Topology (especially differential
structures on spheres), Theory of Modular Forms (especially Eisenstein
series), $p$-Adic Analytic Number Theory (especially $p$-adic $L$-functions)
and Quantum Physics (especially quantum groups). Investigations involving
the Genocchi polynomials and their associated combinatorial relations have
received considerable attention in recent years (see, for details, \cite%
{Araci6}, \cite{Araci2}, \cite{Araci 3}, \cite{Araci 4}, \cite{Kim 3}, \cite%
{Rim}, \cite{Srivastava2} and \cite{Srivastava3}).

Araci \textit{et al.} \cite{Araci2} studied the modified $q$-Genocchi
polynomials which are given by the following generating function: 
\begin{equation}
F_{q}\left( x,t\right) =\sum_{n=0}^{\infty }\mathcal{G}_{n,q}\left( x\right) 
\frac{t^{n}}{n!}=\left[ 2\right] _{q}t\sum_{m=0}^{\infty }\left( -q\right)
^{m}e^{\left( x+\left[ m\right] _{q}\right) t}\text{,}  \label{e3}
\end{equation}%
where the $q$-number $[\lambda ]_{q}$ is given by 
\begin{equation}
\left[ \lambda \right] _{q}:=\frac{1-q^{\lambda }}{1-q}\qquad
(0<q<1;\;\lambda \in \mathbb{C})\text{,}  \label{q-number}
\end{equation}%
so that, obviously, we have 
\begin{equation*}
\lim_{q\rightarrow 1^{-}}\{\left[ \lambda \right] _{q}\}=\lambda \qquad
(\lambda \in \mathbb{C})\text{.}
\end{equation*}%
In the case when $x=0$ in (\ref{e3}), it leads to 
\begin{equation*}
\mathcal{G}_{n,q}\left( 0\right) :=\mathcal{G}_{n,q}\text{,}
\end{equation*}%
that is, to the modified $q$-Genocchi numbers $\mathcal{G}_{n,q}$. In
addition to this, by letting $q\rightarrow 1^{-}$, $G_{n,q}$ reduces to the
Genocchi numbers $G_{n}$: 
\begin{equation*}
\lim_{q\rightarrow 1^{-}}\{\mathcal{G}_{n,q}\}=G_{n}\text{.}
\end{equation*}

The Genocchi numbers $G_{n}(x)$ possess a number of important properties and
are well known in Number Theory. In fact, these numbers are related to the
values at negative integers of the Euler Zeta function defined by (see \cite%
{Kim 10}, \cite{ryoo1}, \cite{Jang}, \cite{Srivastava4}, \cite%
{SrivastavaChoi2012}; see also \cite{Srivastava1}) 
\begin{equation}
\zeta \left( s,x\right) =\sum_{n=0}^{\infty }\frac{\left( -1\right) ^{n}}{%
\left( x+n\right) ^{s}}=\Phi (-1,s,x)  \label{e4}
\end{equation}%
\begin{equation*}
(s\in \mathbb{C};\;x\in \mathbb{C}\setminus \mathbb{Z}_{0}^{-};\;\mathbb{Z}%
_{0}^{-}:=\{0,-1,-2,\cdots \})\text{,}
\end{equation*}%
where $\Phi (z,s,a)$ denotes the widely- and extensively-studied general
Hurwitz-Lerch Zeta function defined by (see, for example, \cite[p. 121 \emph{%
et seq.}]{Srivastava4} and \cite[p. 194 \textit{et seq.}]{SrivastavaChoi2012}%
; see also \cite{HMS2014}, \cite{Srivastava1} and \cite{SSPS}) 
\begin{equation}
\Phi (z,s,a):=\sum_{n=0}^{\infty }\frac{z^{n}}{(n+a)^{s}}  \label{Phi}
\end{equation}%
\begin{equation*}
\big(a\in \mathbb{C}\setminus \mathbb{Z}_{0}^{-};\;s\in \mathbb{C}\quad 
\text{when}\quad |z|<1;\;\Re (s)>1\quad \text{when}\quad |z|=1\big)\text{.}
\end{equation*}

Recently, Kim \cite{Kim 10} defined the $q$-Euler Zeta function as follows: 
\begin{equation}
\zeta _{q}\left( s,x\right) =\left[ 2\right] _{q}\sum_{n=0}^{\infty }\frac{%
\left( -1\right) ^{n}q^{n}}{\left[ x+n\right] _{q}^{s}}\qquad (s\in \mathbb{C%
};\;x\in \mathbb{C}\setminus \mathbb{Z}_{0}^{-})\text{.}  \label{e5}
\end{equation}%
On the other hand, Araci \textit{et al} \cite{Araci2}, Araci introduced the $%
q$-Zeta type function $\widetilde{\zeta }_{q}\left( s,x\right) $ which is
slightly different from Kim's $q$-Zeta function $\zeta _{q}\left( s,x\right) 
$ defined by \eqref{e5}: 
\begin{align}
\widetilde{\zeta }_{q}\left( s,x\right) & :=\frac{1}{\Gamma \left( s\right) }%
\int_{0}^{\infty }t^{s-2}\left\{ -F_{q}\left( x,-t\right) \right\} \mathrm{d}%
t  \notag  \label{e6} \\
& \;=\left[ 2\right] _{q}\sum_{n=0}^{\infty }\frac{\left( -1\right) ^{n}q^{n}%
}{\left( x+\left[ n\right] _{q}\right) ^{s}}\qquad \big(s\in \mathbb{C}%
;\;x\neq -[n]_{q}\;\;\;(n\in \mathbb{N}_{0})\big)\text{,}
\end{align}%
where $F_{q}\left( x,-t\right) $ is given by \eqref{e3}. From (\ref{e3}) and
(\ref{e6}), we find that (see \cite{Araci2}) 
\begin{equation}
\widetilde{\zeta }_{q}\left( -n,x\right) =\frac{\mathcal{G}_{n+1,q}\left(
x\right) }{n+1}\qquad (n\in \mathbb{N}_{0})\text{.}  \label{e7}
\end{equation}%
Moreover, by using (\ref{e5}) and (\ref{e6}), we have 
\begin{equation}
q^{-sx}\widetilde{\zeta }_{q}\left( s,q^{-1}\left[ x\right] _{q^{-1}}\right)
=\zeta _{q}\left( s,x\right) \text{.}  \label{e8}
\end{equation}

The Zeta functions play a crucially important r\^{o}le in Analytic Number
Theory and have applications in such areas as (for example) physics,
probability theory, applied statistics, complex analysis, mathematical
physics, $p$-adic analysis and other related areas. In particular, the Zeta
functions occur within the concept of knot theory, quantum field theory,
applied analysis and number theory (see \cite{Bagdasaryan1}, \cite%
{Bagdasaryan2}, \cite{Cetin}, \cite{Kim 10}, \cite{Kim 6}, \cite{ryoo1}, 
\cite{Jang}, \cite{Srivastava4} and \cite{Srivastava1}).

The distribution formula for the modified $q$-Genocchi polynomials is given
by (see \cite{Araci2}) 
\begin{equation}
\mathcal{G}_{n,q}\left( q^{a}\left[ d\right] _{q}x\right) :=\frac{\left[ d%
\right] _{q}^{n-1}}{\left[ d\right] _{-q}}\sum_{a=0}^{d-1}\left( -1\right)
^{a}q^{a\left( n+1\right) }\mathcal{G}_{n,q^{d}}\left( x+\frac{\left[ a%
\right] _{q}}{q^{a}\left[ d\right] _{q}}\right) \text{, for }d\equiv 1\left( 
\func{mod}2\right) \text{.}  \label{e9}
\end{equation}

Araci \textit{et al.} \cite{Araci 4} derived several new identities for the $%
(h,q)$-Genocchi polynomials and gave symmetric identities of the $(h,q)$%
-Zeta type functions. Yuan He \cite{He} gave symmetric identities for
Carlitz's $q$-Bernoulli numbers (see also \cite{Choi1} and \cite{Choi2}).
Kim also obtained symmetric identities for the $q$-Euler polynomials and
derived the symmetric identities for the $q$-Euler Zeta function (see \cite%
{Kim7}). Simsek \cite{Simsek} gave the complete sum of products of $(h,q)$%
-extension of the Euler polynomials. Bagdasaryan investigated the elementary
evaluation of the Zeta function and presented a real analytic approach to
the values of the Riemann Zeta function (see, for details, \cite%
{Bagdasaryan1} and \cite{Bagdasaryan2}).

The symmetric identity of the Genocchi polynomials is given by Theorem 1
below (see \cite{Cetin}).

\begin{theorem}
\label{thm0} Let $a$ and $b$ be odd integers. Then we have 
\begin{align}
& \sum_{i=0}^{m}\binom{m}{i}a^{i-1}b^{m-i}G_{i}\left( bx\right) S_{m-i}(a) 
\notag  \label{e10} \\
& \qquad \quad =\sum_{i=0}^{m}\binom{m}{i}b^{i-1}a^{m-i}G_{i}\left(
ax\right) S_{m-i}(b)\text{,}
\end{align}%
where 
\begin{equation}
S_{m}(a):=\sum_{j=0}^{a-1}(-1)^{j}\;j^{m}\text{.}  \label{S}
\end{equation}
\end{theorem}

Motivated essentially by some of the aforecited investigations, the
fundamental aim of this paper is to generalize Theorem \ref{thm0} by
presenting an interesting and potentially useful extension of the symmetry
identity (\ref{e10}) to hold true for the modified $q$-Genocchi polynomials
arising from the above-mentioned $q$-Zeta type functions. Several other
related results are also considered.

\section{\textbf{The $q$-Zeta Type Functions}}

In this section, we recall from \eqref{e6} that 
\begin{equation}
\widetilde{\zeta }_{q}\left( s,x\right) =\left[ 2\right] _{q}\sum_{m=0}^{%
\infty }\frac{\left( -1\right) ^{m}q^{m}}{\left( x+\left[ m\right]
_{q}\right) ^{s}}  \label{equation 6}
\end{equation}

In view of (\ref{e8}), we consider (\ref{equation 6}) in the following form: 
\begin{equation}
q^{-asbx-sbj}\widetilde{\zeta }_{q^{a}}\left( s,q^{-a}\left[ bx+\frac{bj}{a}%
\right] _{q^{-a}}\right) =\left[ 2\right] _{q^{a}}\sum_{m=0}^{\infty }\frac{%
\left( -1\right) ^{m}q^{ma}}{\left[ m+bx+\frac{bj}{a}\right] _{q^{a}}^{s}}%
\text{.}
\end{equation}

For non-negative integers $k$ and $i$ such that $m=bk+i$ with $0\leqq i\leqq
b-1$, if we suppose that $a\equiv 1(\func{mod}\;2)$ and $b\equiv 1(\func{mod}%
\;2)$, then we have 
\begin{align}
q^{-asbx-sbj}\widetilde{\zeta }_{q^{a}}\left( s,q^{-a}\left[ bx+\frac{bj}{a}%
\right] _{q^{-a}}\right) & =\left[ a\right] _{q}^{s}\left[ 2\right]
_{q^{a}}\sum_{m=0}^{\infty }\frac{\left( -1\right) ^{m}q^{ma}}{\left[
ma+abx+bj\right] _{q^{a}}^{s}}  \notag  \label{equation 7} \\
& =\left[ a\right] _{q}^{s}\left[ 2\right] _{q^{a}}\sum_{m=0}^{\infty
}\sum_{i=0}^{b-1}\frac{\left( -1\right) ^{i+mb}q^{\left( i+mb\right) a}}{%
\left[ \left( i+mb\right) a+abx+bj\right] _{q^{a}}^{s}}  \notag \\
& =\left[ a\right] _{q}^{s}\left[ 2\right] _{q^{a}}\sum_{i=0}^{b-1}\left(
-1\right) ^{i}q^{ia}\sum_{m=0}^{\infty }\frac{\left( -1\right) ^{m}q^{mba}}{%
\left[ ab\left( m+x\right) +ai+bj\right] _{q^{a}}^{s}}\text{,}
\end{align}%
which readily yields 
\begin{align}
& \sum_{j=0}^{a-1}\left( -1\right) ^{j}q^{jb}q^{-asbx-sbj}\widetilde{\zeta }%
_{q^{a}}\left( s,q^{-a}\left[ bx+\frac{bj}{a}\right] _{q^{-a}}\right)  \notag
\label{equation 8} \\
& \qquad \quad =\left[ a\right] _{q}^{s}\left[ 2\right] _{q^{a}}%
\sum_{j=0}^{a-1}\left( -1\right) ^{j}q^{jb}\sum_{i=0}^{b-1}\left( -1\right)
^{i}q^{ia}\sum_{m=0}^{\infty }\frac{\left( -1\right) ^{m}q^{mba}}{\left[
ab\left( m+x\right) +ai+bj\right] _{q^{a}}^{s}}.
\end{align}

Upon replacing $a$ by $b$ and $j$ by $i$ in (\ref{equation 7}), we get 
\begin{align}
& q^{-asbx-as}\widetilde{\zeta }_{q^{b}}\left( s,q^{-b}\left[ ax+\frac{ai}{b}%
\right] _{q^{-b}}\right)  \notag  \label{222} \\
& \qquad \quad =\left[ b\right] _{q}^{s}\left[ 2\right] _{q^{b}}%
\sum_{j=0}^{a-1}\left( -1\right) ^{j}q^{jb}\sum_{m=0}^{\infty }\frac{\left(
-1\right) ^{m}q^{mba}}{\left[ ab\left( m+x\right) +ai+bj\right] _{q^{a}}^{s}}%
.
\end{align}%
Thus, by applying (\ref{equation 8}) in \eqref{222}, we obtain the following
theorem.

\begin{theorem}
\label{thm1} For any odd integers $a$ and $b$, we have 
\begin{align}
& \frac{\left[ 2\right] _{q^{b}}}{\left[ a\right] _{q}^{s}}%
\sum_{i=0}^{a-1}\left( -1\right) ^{i}q^{ib\left( 1-s\right) }\widetilde{%
\zeta }_{q^{a}}\left( s,q^{-a}\left[ bx+\frac{bi}{a}\right] _{q^{-a}}\right)
\notag \\
& \qquad \quad =\frac{\left[ 2\right] _{q^{a}}}{\left[ b\right] _{q}^{s}}%
\sum_{i=0}^{b-1}\left( -1\right) ^{i}q^{ia\left( 1-s\right) }\widetilde{%
\zeta }_{q^{b}}\left( s,q^{-b}\left[ ax+\frac{ai}{b}\right] _{q^{-b}}\right) 
\text{.}
\end{align}
\end{theorem}

\noindent \textbf{Remark 1.} Upon setting $b=1$ in Theorem \ref{thm1}, we
easily deduce that 
\begin{equation}
\widetilde{\zeta }_{q}\left( s,q^{-1}\left[ ax\right] _{q^{-1}}\right) =%
\frac{\left[ 2\right] _{q}}{\left[ 2\right] _{q^{a}}\left[ a\right] _{q}^{s}}%
\sum_{i=0}^{a-1}\left( -1\right) ^{i}q^{i\left( 1-s\right) }\widetilde{\zeta 
}_{q^{a}}\left( s,q^{-a}\left[ x+\frac{i}{a}\right] _{q^{-a}}\right) \text{.}
\label{equation 9}
\end{equation}

Taking $a=2$ in (\ref{equation 9}), we derive the following Corollary.

\begin{corollary}
For any odd integer $a$, we have 
\begin{equation}
\widetilde{\zeta }_{q}\left( s,q^{-1}\left[ 2x\right] _{q^{-1}}\right) =%
\frac{\left[ 2\right] _{q}}{\left[ 2\right] _{q^{2}}\left[ 2\right] _{q}^{s}}%
\left[ \widetilde{\zeta }_{q^{2}}\left( s,q^{-2}\left[ x\right]
_{q^{-2}}\right) -q^{b\left( 1-s\right) }\widetilde{\zeta }_{q^{2}}\left(
s,q^{-2}\left[ x+\frac{1}{2}\right] _{q^{-2}}\right) \right] \text{.}
\end{equation}
\end{corollary}

\noindent \textbf{Remark 2.} If we take $s=-n$ in Theorem \ref{thm1}, we get
the following symmetric property of the modified $q$-Genocchi polynomials.

\begin{theorem}
\label{thm2} For any odd integers $a$ and $b$, we have 
\begin{align}
& \left[ 2\right] _{q^{b}}\left[ a\right] _{q}^{n-1}\sum_{i=0}^{a-1}\left(
-1\right) ^{i}q^{ib\left( n+1\right) }\mathcal{G}_{n,q^{a}}\left( q^{-a}%
\left[ bx+\frac{bi}{a}\right] _{q^{-a}}\right) \\
& \qquad \quad =\left[ 2\right] _{q^{a}}\left[ b\right] _{q}^{n-1}%
\sum_{i=0}^{b-1}\left( -1\right) ^{i}q^{ia\left( n+1\right) }\mathcal{G}%
_{n,q^{b}}\left( q^{-b}\left[ ax+\frac{ai}{b}\right] _{q^{-b}}\right) \text{.%
}
\end{align}
\end{theorem}

We now take $b=1$ and replace $x$ by $\frac{x}{a}$ in Theorem \ref{thm2}. We
thus restate the distribution formula for the modified $q$-Genocchi
polynomials as follows: 
\begin{equation}
\mathcal{G}_{n,q}\left( -\left[ -x\right] _{q}\right) =\frac{\left[ 2\right]
_{q}}{\left[ 2\right] _{q^{a}}}\left[ a\right] _{q}^{n-1}\sum_{i=0}^{a-1}%
\left( -1\right) ^{i}q^{i\left( n+1\right) }\mathcal{G}_{n,q^{a}}\left(
q^{-a}\left[ \frac{x+i}{a}\right] _{q^{-a}}\right) \text{ }\left( 2\nmid
a\right) \text{.}
\end{equation}

We next find from (\ref{e3}) that 
\begin{align*}
\sum_{n=0}^{\infty }\mathcal{G}_{n,q}\left( x+y\right) \frac{t^{n}}{n!}& =%
\left[ 2\right] _{q}t\sum_{m=0}^{\infty }\left( -q\right) ^{m}e^{\left( x+y+%
\left[ m\right] _{q}\right) t} \\
& =\left( \sum_{m=0}^{\infty }y^{m}\frac{t^{m}}{m!}\right) \left(
\sum_{n=0}^{\infty }\mathcal{G}_{n,q}\left( x\right) \frac{t^{n}}{n!}\right) 
\text{,}
\end{align*}%
which, by applying the Cauchy product, yields 
\begin{equation}
\sum_{n=0}^{\infty }\mathcal{G}_{n,q}\left( x+y\right) \frac{t^{n}}{n!}%
=\sum_{n=0}^{\infty }\left( \sum_{k=0}^{n}\binom{n}{k}\mathcal{G}%
_{k,q}\left( x\right) y^{n-k}\right) \frac{t^{n}}{n!}\text{.}  \label{111}
\end{equation}%
Thus, by comparing the coefficients of $\frac{t^{n}}{n!}$ on both sides of
this last equation \eqref{111}, we get the following Corollary.

\begin{corollary}
For $n\in \mathbb{N}_{0}$, we obtain 
\begin{equation}
\mathcal{G}_{n,q}\left( x+y\right) =\sum_{k=0}^{n}\binom{n}{k}\mathcal{G}%
_{k,q}\left( x\right) y^{n-k}\text{.}  \label{equation 10}
\end{equation}
\end{corollary}

By using Theorem \ref{thm2} and (\ref{equation 10}), we can derive Theorem %
\ref{Theorem 2.5} below.

\begin{theorem}
\label{Theorem 2.5} For any odd integers $a$ and $b$, we have 
\begin{align}
& \left[ 2\right] _{q^{b}}\left[ a\right] _{q}^{n-1}\sum_{k=0}^{n}\binom{n}{k%
}\left[ a\right] _{q^{-1}}^{k-n}\left[ b\right] _{q^{-1}}^{n-k}\mathcal{G}%
_{k,q}\left( q^{-a}\left[ bx\right] _{q^{-a}}\right) S_{n-k:q^{-b}}^{\left(
n+1\right) }\left( a\right)   \notag \\
& \qquad =\left[ 2\right] _{q^{a}}\left[ b\right] _{q}^{n-1}\sum_{k=0}^{n}%
\binom{n}{k}\left[ b\right] _{q^{-1}}^{k-n}\left[ a\right] _{q^{-1}}^{n-k}%
\mathcal{G}_{k,q}\left( q^{-b}\left[ ax\right] _{q^{-b}}\right)
S_{n-k:q^{-a}}^{\left( n+1\right) }\left( b\right) 
\end{align}%
where 
\begin{equation}
S_{m:q}^{\left( j\right) }\left( a\right) :=\sum_{i=0}^{a-1}\left( -1\right)
^{i}q^{ji}\left[ i\right] _{q}^{m}\text{.}
\end{equation}
\end{theorem}

\noindent \textbf{Remark 3.} Letting $q\rightarrow 1^{-}$ in Theorem \ref%
{Theorem 2.5}, we can deduce the known symmetry identity (\ref{e10}).

\section{\textbf{Concluding Remarks and Observations}}

In this article, we have derived several symmetric properties of the $q$%
-Zeta type function $\widetilde{\zeta }_{q}\left( s,x\right) $ defined by %
\eqref{e6}. As applications of these properties, we give new interesting
symmetry identities for the modified $q$-Genocchi polynomials $\mathcal{G}%
_{n,q}\left( x\right) $ which are defined by \eqref{e3}. In the limit when $%
q\rightarrow 1^{-}$, this last result (Theorem 4) is shown to yield the
known symmetry identity \eqref{e10} for the Genocchi polynomials $G_{n}(x)$.

\end{document}